\documentclass{article}
\usepackage{amssymb}
\usepackage[centertags]{amsmath}

\begin{document}

\title{The renormalization method from continuous to discrete dynamical systems: asymptotic solutions, reductions and invariant manifolds}
\author{ Cheng-shi Liu \\Department of Mathematics\\Northeast Petroleum University\\Daqing 163318, China
\\Email: chengshiliu-68@126.com}

 \maketitle

\begin{abstract}

 The renormalization method based on the Taylor expansion for asymptotic analysis of differential equations is generalized to difference equations. The proposed renormalization method is based on the Newton-Maclaurin expansion. Several basic theorems on the renormalization method are proven. Some interesting applications are given, including asymptotic solutions of quantum anharmonic oscillator and discrete boundary layer, the reductions and invariant manifolds of some discrete dynamics systems. Furthermore, the homotopy renormalization method based on the Newton-Maclaurin expansion is proposed and applied to those difference equations including no a small parameter.


\end{abstract}

\section{Introduction}

 The renormalization method based on the Taylor expansion (TR, for simplicity) and corresponding homotopy renormalization method (HTR for simplicity) have been obtained and applied to a large number of the perturbation and non-perturbation differential equations to give their global valid asymptotic solutions[1].  In particular, Goldenfeld et al's renormalization group (RG, for simplicity) method [2-4] and its geometrical formulations[5-7] can be derived from our renormalization method.  The RG method unifies
some perturbation methods [8-13] including singular and reductive
perturbation theory and hence shows its advantages, and  have been
developed and applied to many differential equations[2-7, 14-21].

Naturally, we should consider the renormalization theory of difference equations. Difference equations[22-25] and its perturbation theory had been studied extensively[26-38]. For example, Huston[29] used the Krylov-Bogoljibov method to difference equations, Marathe and Chaterjee[36] used harmonic balance and multiple scales methods to nonlinear periodic structure, Horssen et al[37,38] proposed a modified multiple scales method to solve the difference equations, and so forth. In particular, Kunihiro and Matsukidaira [39] generalized the RG method to the discrete systems. However, both of the standard RG method and its geometrical interpretation are still not
satisfied since their mathematical foundations are not very clear.

In the paper,  based on the Newton-Maclaurin expansion, we obtain the renormalization method and homotopy renormalization method for difference equations, and give some  applications, including the global valid asymptotic solutions for some  regular and singular perturbation problems, reductions and invariant manifolds of discrete systems, and some difference equations without a small parameter. Comparing with those routine methods such as multiple scales method and RG method, our renormalization method is clear in theory and more simple in practice.
  The biggest advantage of
our method is that the secular terms can be automatically
eliminated, and hence it does not
require the asymptotic matching and not need the introduction of
auxiliary renormalization parameters.

This paper is organized as follows. In section 2, we summarize some main results on renormalization method based on the Taylor expansion and give several further new applications such as quantum anharmonic oscillator. In section 3,  we obtain the renormalization method based on the Newton-Maclaurin expansion for asymptotic
analysis of difference equations.  In section 4, we use the proposed renormalization method
 to asymptotic analysis for some  difference equations such as discrete boundary layer problems, the invariant manifold and reduction equations of discrete dynamics systems. In section 5, we propose the homotopy renormalization method and solve some difference equations. The last section is a short conclusion.

\section{The renormalization method based on the Taylor expansion and some further applications}

In the section, we summarize some main results on the renormalization method based on the Taylor expansion(see [1] for more details and applications) and give several new applications. For a differential
equation
\begin{equation}
N(y)=\epsilon M(y),
\end{equation}
where $N$ and $M$ are in general linear or nonlinear operators.
Assuming that the solution can be expanded as a power series of the
small parameter $\epsilon$
\begin{equation}
y=y_0+y_1\epsilon+\cdots+y_n\epsilon^n+\cdots,
\end{equation}
and substituting it into the above equation yields the equations of
$y_n$'s such as
\begin{equation}
N(y_0)=0,
\end{equation}
and
\begin{equation}
N_1(y_1)=M_1(y_0),\cdots,N_k(y_k)=M_k(y_{k-1}),
\end{equation}
for some operator $N_k$ and $M_k$. By the first equation, we give
the general solution of $y_0$ including some integral constants $A$
and $B$ and so on. In general the number of the integral constants
is equal to the order of the differential equation. Then we find the
special solutions of the $y_k$ which sometimes include also some
integral constants, and expand them as a power series at a general
point $t_0$ [41-44]. Therefore, by rearranging the summation of these
series, we obtain the final solution
\begin{equation}
y(t,t_0)=\sum_{n=0}^{+\infty}Y_n(t_0,\epsilon)(t-t_0)^n.
\end{equation}
This is the most important formula for our theory from which we can
give every thing of the standard RG method and more. In fact, the main results on our renormalization method in [1] can be we summarized as follows:

 \textbf{Theorem 2.1}. The exact solution of the
Eq.(1) is just
\begin{equation}
y(t)=Y_0(t,\epsilon),
\end{equation}
and furthermore, we have
\begin{equation}
Y'_{n-1}(t,\epsilon)=nY_n(t,\epsilon),n=1,2,\cdots.
\end{equation}

This theorem plays a key role in our theory. The first
formula (6) tell us that we only need to find $Y_0$,  and the second formula (7)
tell us how to determine $Y_0$. In fact, if there exist $m$ integral
constants in $Y_0$ to be determined, we will need the second formula (7)
to provide $m$ natural relations to solve these constants. In
principle, we need $m$ equations to determine $m$ unknown integral
constants.

\textbf{Theorem 2.2}. From the expansion (5), the standard RG equation(see [3])
holds, that is
\begin{equation}
\frac{\partial}{\partial t_0}y(t,t_0)=0.
\end{equation}

\textbf{Definition 1}: If there exist $m$ unknown constants, the
equations
\begin{equation}
Y'_{n-1}(t,\epsilon)=nY_n(t,\epsilon),n=1,2,\cdots,m.
\end{equation}
are called renormalization equations by which the $m$ constants can be
determined.

\textbf{Remark 1}. In general, the renormalization equations are so complicated
that we need some reduced treatments. In practice, we in general
only need the first equation
\begin{equation}
Y'_0(t,\epsilon)=Y_1(t,\epsilon),
\end{equation}
and split it to find each integral constant.

Next we derive out the assumption $t_0=t$ in Kunihiro's geometrical
formulation[5-7]. Considering the solution (5) as a local solution,
the global solution is the envelop of these local solutions by
eliminating the parameter $t_0$ from  equation (5) and the RG
equation $\frac{\partial}{\partial t_0}y(t,t_0)=0$. However, the RG
equation is in general so complicated that Kunihiro[5-7] has to add an
auxiliary assumption $t_0=t$. In the following we prove that this is
a fact but not an assumption.

 \textbf{Theorem 2.3}. As an envelop of the local solutions, the global solution is tangent with the local
 solution at the point $t_0=t$.

 As pointed out in [1], we call this renormalization method based on the Taylor expansion the Taylor renormalization method, for simplicity, TR method. A further development of TR method is called the homotopy renormalization method (HTR method for simplicity). In [1], a large number of applications of TR and HTR methods have been given. In the following we give further new applications of TR method.

\textbf{Example 1}. Consider the perturbation equation([3]) $y''(t)+y(t)-\epsilon t y(t)=0$.

By the renormalization method, taking
\begin{equation}
y=\sum_{n=0}^{+\infty}y_n\epsilon^n,
\end{equation}
and substituting it into the above equation yields
\begin{equation}
y=R\cos(t+\theta)+\epsilon\{\frac{R}{2}(1+t_0)-\frac{R}{4}\frac{t_0}{\sin^2(t_0+\theta)}+\cdots\}(t-t_0)
+O(\epsilon^2),
\end{equation}
which includes secular terms. By considering $R$ and $\theta$ as the functions of $t_0$, the renormalization equation
\begin{equation*}
\frac{\partial}{\partial
t_0}(R\cos(t_0+\theta))
\end{equation*}
\begin{equation}
=\epsilon\{\frac{R}{2}(1+t_0)-\frac{R}{4}\frac{t_0}{\sin^2(t_0+\theta)}+\cdots\},
\end{equation}
gives a closed equations system
\begin{equation}
\theta'=-\frac{\epsilon}{2} (1+t_0),
\end{equation}
\begin{equation}
R'=0.
\end{equation}
Solving the above equations give the solutions
\begin{equation}
R(t_0)=R_0,
\end{equation}
\begin{equation}
\theta(t_0)=-\frac{\epsilon}{4}(1+t_0)^2+\theta_0,
\end{equation}
where $R_0$ and $\theta_0$ are integral constants. Therefore, the global asymptotic approximate solution is given by
\begin{equation}
y(t)=R_0\cos(t-\frac{\epsilon}{4}(1+t)^2+\theta_0)+O(\epsilon).
\end{equation}

\textbf{Example 2}. Consider the equation ([3]) $\epsilon y''(t)+2y'(t)+\exp(y(t))=0$.

Taking the transformations $t=\epsilon x, y(t)=z(x)$, we have
\begin{equation}
 z''(x)+2z'(x)+\epsilon \exp (z(x))=0.
\end{equation}
By taking
\begin{equation}
z=\sum_{n=0}^{+\infty}z_n\epsilon^n,
\end{equation}
and substituting it into the above equation yields
\begin{equation}
z''_0+2z'_0=0,
\end{equation}
and
\begin{equation}
z''_1+2z'_1=- \exp (z_0(x)),
\end{equation}
and so on. The general solution of the first equation is
\begin{equation}
z_0=A+B\exp(-2x).
\end{equation}
The general solution of the above equation can be also easily
obtained. But here we find its power series solution. Assume that
\begin{equation}
z_1=h_0+h_1(x-x_0)+h_2(x-x_0)^2+h_3(x-x_0)^3+\cdots
\end{equation}
and expand $\exp(B\exp(-2x))$ as
\begin{equation}
\exp(B\exp(-2x))=\exp(B\exp(-2x_0))-2B\exp(-2x_0)\exp(B\exp(-2x_0))(x-x_0)+\cdots.
\end{equation}
Substituting (28) into Eq.(25) yields
\begin{equation}
2h_1+2h_2=-\exp(A)\exp(B\exp(-2x_0)),
\end{equation}
\begin{equation}
4h_2+6h_3=2B\exp(A)\exp(-2x_0)\exp(B\exp(-2x_0)).
\end{equation}
Letting $h_2=0$ gives
\begin{equation}
h_1=-\frac{1}{2}\exp(A)\exp(B\exp(-2x_0)),h_3=\frac{1}{3}B\exp(A)\exp(-2x_0)\exp(B\exp(-2x_0)).
\end{equation}
Furthermore, by taking $h_0=0$, we have
\begin{equation*}
z_1=-\frac{1}{2}\exp(A)\exp(B\exp(-2x_0))(x-x_0)
\end{equation*}
\begin{equation}
+\frac{1}{3}B\exp(A)\exp(-2x_0)\exp(B\exp(-2x_0))(x-x_0)^3+\cdots,
\end{equation}
and hence
\begin{equation*}
z=A+B\exp(-2x)-\frac{\epsilon}{2}\exp(A)\{1+B\exp(-2x_0)+\cdots\}(x-x_0)
\end{equation*}
\begin{equation}
+\frac{\epsilon}{3}B\exp(A)\exp(-2x_0)\exp(B\exp(-2x_0))(x-x_0)^3+\cdots.
\end{equation}
The renormalization equation is
\begin{equation}
(A+B\exp(-2x_0))'=-2
B\exp(-2x_0)-\frac{\epsilon}{2}\exp(A)(1+B\exp(-2x_0)),
\end{equation}
where we have used the Taylor expansion formula
\begin{equation}
\exp(-2x)=\exp(-2x_0)-2\exp(-2x_0)(x-x_0)+\cdots.
\end{equation}
According to the renormalization equation, we can take a closed equations system
\begin{equation}
A'=-\frac{\epsilon}{2}\exp(A), B'=-\frac{\epsilon}{2}\exp(A)B,
\end{equation}
whose solutions are
\begin{equation}
A(x)=\ln\frac{2}{\epsilon x+C}, B(x)=\frac{B_0}{\epsilon x+C}.
\end{equation}
Therefore, the global approximate solution is given by
\begin{equation}
z(x)=\ln\frac{2}{\epsilon x+C}+\frac{B_0}{\epsilon x+C}\exp(-2x).
\end{equation}
Respectively, we have the global asymptotic approximate solution
\begin{equation}
y(t)=\ln\frac{2}{t+C}+\frac{B_0}{ t+C}\exp(-2t/\epsilon)+O(\epsilon).
\end{equation}
According to the boundary conditions, we can determine the values of
$C$ and $B_0$.

\textbf{Example 3}. Consider the quantum anharmonic oscillator
\begin{equation}
\phi''(x)-\frac{1}{4}x^2\phi(x)-\frac{1}{4}\epsilon x^4\phi(x)=-E\phi(x),
\end{equation}
with the boundary conditions
\begin{equation}
\phi(\pm\infty)=0.
\end{equation}

By the WKB method, it follows that for large $x$,
\begin{equation}
\phi(x)\sim \exp(-\sqrt\epsilon |x|^3/6).
\end{equation}
Bender and Betencourt [40] used the Rayleigh-Schr\"{o}dinger perturbation theory to obtain the result. Kunihiro [21] reconstructed this result by the geometrical formula of the renormalization group method. Here we use the TR method to deal with it.

By the Rayleigh-Schr\"{o}dinger perturbation method, we represent both the eigenfunction and eigenvalue as the power series in $\epsilon$
 \begin{equation}
\phi(x)=\sum_{k=0}^{+\infty}\phi_k(x)\epsilon^k,
\end{equation}
 \begin{equation}
E=\sum_{k=0}^{+\infty}E_k\epsilon^k.
\end{equation}
Substituting them into the Schr\"{o}dinger equation (40) yields
\begin{equation}
\phi_0''(x)-\frac{1}{4}x^2\phi(x)=-E_0\phi_0(x),
\end{equation}
\begin{equation}
\phi_1''(x)-\frac{1}{4}x^2\phi_1(x)-\frac{1}{4} x^4\phi_0(x)=-E_1\phi_0(x)-E_0\phi_1(x),
\end{equation}
and so on. Solving the first equation gives
\begin{equation}
\phi_0(x)=A\mathrm{e}^{-x^2/4},
\end{equation}
\begin{equation}
E_0=\frac{1}{2}.
\end{equation}
Further, we take for $k\geq1$
\begin{equation}
\phi_k(x)=A\mathrm{e}^{-x^2/4}P_k(x),
\end{equation}
and hence
\begin{equation}
P_k''(x)-xP_k'(x)=\frac{x^4}{4}P_{k-1}(x)-\sum_{j=1}^{k-1}P_j(x)E_{k-j},
\end{equation}
\begin{equation}
E_k=\frac{1}{2\sqrt\pi}\int_{-\infty}^{+\infty}\mathrm{e}^{-x^2/4}\{\frac{x^4}{4}P_{k-1}(x)
-\sum_{j=1}^{k-1}P_j(x)E_{k-j}\}\mathrm{d}x,
\end{equation}
where $P_0(x)=1$. For $k=1$, we give the solution at a general point $x_0$
\begin{equation}
\phi_1(x)=A\mathrm{e}^{-x^2/4}\{P_1(x)-P_1(x_0)\},
\end{equation}
where
\begin{equation}
P_1(x)=-\frac{1}{16}x^4-\frac{3}{8}x^2.
\end{equation}
Similarly, for $k=2$, we have
\begin{equation}
\phi_2(x)=A\mathrm{e}^{-x^2/4}\{P_2(x)-P_2(x_0)-P_1(x_0)(P_1(x)-P_1(x_0))\},
\end{equation}
where
\begin{equation}
P_2(x)=\frac{1}{512}x^8+\frac{11}{192}x^6+\frac{31}{128}x^4+\frac{1}{16}x^2.
\end{equation}
Therefore, the solution can be given by
\begin{equation}
\phi(x)=A\mathrm{e}^{-x^2/4}\{1+\epsilon(P_1(x)-P_1(x_0))+\epsilon^2(P_2(x)-P_2(x_0)-P_1(x_0)(P_1(x)-P_1(x_0)))\}+\cdots,
\end{equation}
and hence
\begin{equation}
\phi(x)=\mathrm{e}^{-x^2/4}\{A+A(\epsilon P'_1(x_0)+\epsilon^2(P_2'(x_0)-P_1(x_0)P_1'(x_0)))(x-x_0))\}+\cdots.
\end{equation}
By the TR method, we assume that $A$ is a function of $x_0$, then the renormalization equation is
\begin{equation}
A'(x_0)=A(x_0)\{\epsilon P'_1(x_0)+\epsilon^2(P_2'(x_0)-P_1(x_0)P_1'(x_0))+\cdots\},
\end{equation}
by which we get
\begin{equation}
A(x_0)=A_0\exp(\epsilon P'_1(x_0)+\epsilon^2(P_2'(x_0)-P_1(x_0)P_1'(x_0))+\cdots),
\end{equation}
where $A_0$ is an integral constant. Furthermore, according to the TR method, the solutions is just
\begin{equation}
\phi(x)=A(x)\mathrm{e}^{-x^2/4}=A_0\mathrm{e}^{-x^2/4}\exp(\epsilon P'_1(x)+\epsilon^2(P_2'(x)-P_1(x)P_1'(x))+\cdots).
\end{equation}
Using the results of $P_k(x)$ in [39],  we can take the highest order terms in $\phi_k(x)$ for $1\leq k\leq4$, and give
\begin{equation}
\phi(x)=A_0\exp\{-\frac{x^2}{4}(1+2\epsilon x^2+\frac{17}{12}\epsilon^2 x^4+\frac{5}{12}\epsilon^3 x^6
+\frac{47}{1152}\epsilon^4 x^8)+\cdots \}.
\end{equation}
It follows that
\begin{equation}
\phi(x)=A_0\exp\{\frac{1}{6\epsilon}[1-(1+\epsilon x^2)^{3/2}]\}\sim \mathrm{e}^{-\sqrt\epsilon|x|^3/6},
\end{equation}
for $x\rightarrow\pm\infty$.

   By another trick in [40], we take four terms and derive
\begin{equation}
\phi(x)=A_0\exp\{-\frac{x^2}{4}(1+\frac{1}{4}\epsilon x^2-\frac{1}{24}\epsilon^2 x^4+\frac{1}{64}\epsilon^3 x^6
-\frac{1}{128}\epsilon^4 x^8)^{1/8} \}\sim \mathrm{e}^{-\sqrt\epsilon|x|^3/5.96663}.
\end{equation}

The analysis about the solution can be found in [40].

\section{The renormalization method based on the Newton-Maclaurin expansion}
 For a sequence $y(n)$, we can expand it as the Newton-Maclaurin series
\begin{equation}
y(n)=y(m)+\binom{n-m}{1}\Delta y(m)+\binom{n-m}{2}\Delta^2 y(m)+\cdots+\binom{n-m}{r}\Delta^r y(m)+\cdots,
\end{equation}
where we define
\begin{equation}
\Delta y(n)=y(n+1)-y(n),
\end{equation}
\begin{equation}
\Delta^r y(n)=\Delta(\Delta^{r-1}y(n)),  r=1,2,\cdots.
\end{equation}
We will prove the formula (64) in the following theorem 3.1.

Consider a difference equation
\begin{equation}
N(y(n))=\epsilon M(y(n)),
\end{equation}
where $N$ and $M$ are in general linear or nonlinear difference operators.
Assume that the solution can be expanded as a power series of the
small parameter $\epsilon$
\begin{equation}
y(n)=y_0(n)+y_1(n)\epsilon+\cdots+y_k(n)\epsilon^k+\cdots.
\end{equation}
Substituting it into the above equation and equating coefficients of the powers of $\epsilon$ yield
 the equations of
$y_k(n)$'s such as
\begin{equation}
N(y_0(n))=0,
\end{equation}
and
\begin{equation}
N_1(y_1(n))=M_1(y_0(n)),\cdots,N_k(y_k(n))=M_k(y_{k-1}(n)),\cdots,
\end{equation}
for some operators $N_k$ and $M_k$. By the equation (69), we give
the general solution of $y_0$ including some undetermined constants $A$
and $B$ and so on. In general, the number of the undetermined constants
is equal to the order of the difference equation. Next, we find the
particular solutions of $y_k$ which sometimes include also some
undetermined constants, and expand them as the  series at a general
point $m$ as follows,
\begin{equation}
y_j(n)=\sum_{k=0}^{+\infty}\Delta^{k}y_j(m) \binom{n-m}{k},j=0,1,2,\cdots.
\end{equation}
By rearranging the summation of these series, we obtain
the final solution
\begin{equation}
y(n,m)=\sum_{k=0}^{+\infty}Y_k(m,\epsilon) \binom{n-m}{k},
\end{equation}
where
\begin{equation}
Y_k(m,\epsilon)=\sum_{j=0}^{+\infty}\Delta^ky_{j}(m)\epsilon^j,
n=0,1,\cdots.
\end{equation}
The formula (72) is the most basic formula from which we have

(i). $y(m)=y(m,m)=Y_0(m,\epsilon)$;

(ii).$Y_k(m,\epsilon)=\Delta^kY_{0}(m,\epsilon),k=1,2,\cdots$.

These results are the simple facts in the Newton-Maclaurin expansion but it is
the most important result for our theory by the following reasons.
In fact, the formula (i) tells us that the solution is exactly given by the
first term $Y_0$ of the expansion when we consider $m$ as a
parameter, and hence all other terms need not be considered at all!
However, in general, this first term includes some integral
constants which need to be determined.  The point in our theory is to take $r$
relations in case (ii) as the equations to determinate the $r$
unknown parameters. In what following,  we give the proofs of these results.

 \textbf{Theorem 3.1}. The exact solution of the
Eq.(67) is just
\begin{equation}
y(n)=Y_0(n,\epsilon),
\end{equation}
and furthermore, we have
\begin{equation}
\Delta Y_{k-1}(n,\epsilon)=Y_k(n,\epsilon),n=1,2,\cdots.
\end{equation}

\textbf{Proof}. We assume that $y(n)$ can be expanded at a general point $m$ as follows
\begin{equation}
y(n)=a_0+a_1\binom{n-m}{1}+a_2\binom{n-m}{2}+\cdots+a_k\binom{n-m}{k}+\cdots.
\end{equation}
By taking $n=m$, we have $a_0=y(m)$. Next, we take the difference operation in two sides of the above equation and get
\begin{equation}
\Delta y(n)=a_1\binom{n-m}{0}+a_2\binom{n-m}{1}+\cdots+a_k\binom{n-m}{k-1}+\cdots,
\end{equation}
where we use the formula
\begin{equation}
\Delta\binom{n-m}{k}=\binom{n-m}{k-1}.
\end{equation}
Furthermore, we take $n=m$ to give $a_1=\Delta y(m)$. By mathematical induction, we get the formula. From the formula, the conclusion of the theorem is obvious.  We complete the proof.

\textbf{Theorem 3.2}. From the expansion (72), we have
\begin{equation}
\Delta_my(n,m)=0,
\end{equation}
where $\Delta_m$ is the partial difference operator, that is, $\Delta_my(n,m)=y(n,m+1)-y(n,m)$.

\textbf{Proof}. Firstly, we give two useful results which can be derived easily from the definition of $\Delta$ operator:

(i). $\Delta (f(n)g(n))=g(n+1)\Delta f(n)+f(n)\Delta g(n)$;

(ii). $\Delta_m \binom {n-m}{k}=-\binom {n-m-1}{k-1}$, for $k\geq1$.

Then, from the above two results, we have
\begin{equation}
\Delta_m\{\Delta^k y(m)\binom{n-m}{k}\}=\Delta^{k+1} y(m)\binom{n-m-1}{k}-\Delta^k y(m)\binom{n-m-1}{k-1}.
\end{equation}
Furthermore, by direct computation, we give
\begin{equation*}
\Delta_my(n,m)=\Delta y(m)\binom{n-m-1}{0}+\sum_{k=1}^{\infty}\{\Delta^{k+1} y(m)\binom{n-m-1}{k}
\end{equation*}
\begin{equation}
-\Delta^k y(m)\binom{n-m-1}{k-1}\}=0.
\end{equation}
The proof is completed.

\textbf{Definition 2}. We call the relations
\begin{equation}
Y_k(m,\epsilon)=\Delta^kY_{0}(m,\epsilon),k=1,2,\cdots
\end{equation}
as the renormalization equations.

According to the above formulas, we give the main steps of the renormalization method as follows:

Firstly, $Y_0(m,\epsilon)$ is exact solution by expansion (72),
in which there are some constants to be determined.

Secondly, the natural relation $\Delta Y_0=Y_1$ (if there are $r$ unknown
constants, we will take $r$ renormalization equations
$Y_k(m,\epsilon)=\Delta^kY_{0}(m,\epsilon),k=1,2,\cdots,r$)
gives the renormalization equations satisfied by unknown constants in $Y_0$, where these unknown constants
are considered as the functions of variable $m$.

Thirdly, by renormalization equations, we solve out the unknown
constants and substitute them into $Y_0$ to give the asymptotic
solution.

We must emphasize that if there are $r$ unknown integral constants,
we will need $r$ renormalization equations $\Delta Y_{k-1}(m)=Y_k(m)$ to get a
closed equations system. In practice, we usually need only one
renormalization equation, that is, the first renormalization
equation $\Delta Y_{0}(m)=Y_1(m)$ to get the closed equations system by some
approximations or other balance relations.

\textbf{Remark 2}. For the system of finite or infinite-dimensional
ordinary difference equations, we only need to replace the
corresponding scale functions by the vector functions to give the whole
theory. For the partial difference equations, we can easily give the corresponding renormalization method by taking the multi-variables Newton-Maclaurin expansion.

We also need to consider the following asymptotic error property of renormalization solutions. Here we only take the first order vector difference equations as an example to give the theorem.

\textbf{Theorem 3.3}. Consider vector difference equation
\begin{equation}
\Delta Y(n)=F(Y(n), \epsilon),
 \end{equation}
where $Y$ is a vector and $F$ is a vector value function.  We assume that $\widetilde{Y}(n,m)$ is the local solution at the general point $m$ and satisfies
 \begin{equation}
\Delta \widetilde{Y}(n,m)=F(\widetilde{Y}(n,m), \epsilon)+O(\epsilon^k).
 \end{equation}
Then we have
 \begin{equation}
\Delta \widetilde{Y}(m,m)=F(\widetilde{Y}(m,m), \epsilon)+O(\epsilon^k).
 \end{equation}

 \textbf{Proof}. By direct computation, we have
\begin{equation}
\Delta \widetilde{Y}(m,m)=\Delta_n \widetilde{Y}(n,m)|_{n=m}+\Delta_m \widetilde{Y}(n,m)|_{n=m}.
 \end{equation}
Then, according to the theorem 3.2, we have
\begin{equation}
\Delta \widetilde{Y}(m,m)=\Delta_n \widetilde{Y}(n,m)|_{n=m}=F(\widetilde{Y}(m,m), \epsilon)+O(\epsilon^k).
 \end{equation}
The proof is completed.

\textbf{An illustration example}. Consider the equation
 \begin{equation}
 y(n+2)+\epsilon y(n+1)+y(n)=0,
 \end{equation}
where $\epsilon$ is a positive small parameter. Assume that the solution can be expanded as a power series of the
small parameter $\epsilon$
\begin{equation}
y(n)=y_0(n)+y_1(n)\epsilon+\cdots+y_k(n)\epsilon^k+\cdots.
\end{equation}
Substituting it into the above equation yields the equations of
$y_k(n)$'s such as
\begin{equation}
y_0(n+2)+y_0(n)=0,
\end{equation}
and
\begin{equation}
y_1(n+2)+y_1(n)=-y_0(n+1),
\end{equation}
and so forth. Solving the first equation gives
\begin{equation}
y_0(n)=Ai^n+B(-i)^n.
\end{equation}
Then the second equation becomes
\begin{equation}
y_1(n+2)+y_1(n)=-Ai^{n+1}-B(-i)^{n+1}.
\end{equation}
By the variation of constant method, we  obtain
\begin{equation}
y_1(n)=A_1i^{n}+B_1(-i)^{n}+\frac{i}{2}(A-(-1)^mB)i^n(n-m)+\frac{i}{2}(A(-1)^m-B)(-i)^n(n-m).
\end{equation}
Therefore, the solution is given by
\begin{equation*}
y(n)=Ai^n+B(-i)^n+\epsilon A_1i^{n}+\epsilon B_1(-i)^{n}+\frac{i}{2}\epsilon(A-(-1)^mB)i^n(n-m)
\end{equation*}
\begin{equation}
+\frac{i}{2}\epsilon(A(-1)^m-B)(-i)^n(n-m)+O(\epsilon^2),
\end{equation}
which includes secular terms.

By the renormalization method, considering four constants $A$, $A_1$, $B$ and $B_1$ as the functions of variable $m$, and expanding $i^n$ and $(-i)^n$ as the Newton-Maclaurin series at the point $m$, the renormalization equations can be taken as
\begin{equation}
\Delta A(m)=\frac{i}{2}\epsilon A(m), \Delta B(m)=-\frac{i}{2}\epsilon B(m),
\end{equation}
\begin{equation}
\Delta A_1(m)=\Delta B_1(m)=0.
\end{equation}
Solving the renormalization equations, we have
\begin{equation}
A(m)=A_0(1+\frac{i}{2}\epsilon)^m\sim A_0\mathrm{e}^{\frac{m}{2}\epsilon i},
\end{equation}
\begin{equation}
B(m)=B_0(1-\frac{i}{2}\epsilon)^m\sim B_0\mathrm{e}^{-\frac{m}{2}\epsilon i},
\end{equation}
\begin{equation}
A_1=A_{10},B_1=B_{10},
\end{equation}
where $A_0, A_{10},B_0$ and $B_{10}$ are  arbitrary constants. Then, by taking $A_{10}=B_{10}=0$ and $A_0=\overline{B_0}$,  we give the global asymptotic solution
\begin{equation}
y(n)=A_0\mathrm{e}^{\frac{n}{2}(\epsilon+\pi) i}+B_0\mathrm{e}^{-\frac{n}{2}(\epsilon+\pi) i}+O(\epsilon^2)\\
=C_0\cos\frac{n}{2}(\epsilon+\pi)+D_0\sin\frac{n}{2}(\epsilon+\pi)+O(\epsilon^2),
\end{equation}
where $C_0=A_0+\overline{B_0}, D_0=\mathrm{i}(A_0-\overline{B_0})$.

\textbf{Remark 3}. Furthermore, we can also generalize the theory in [41-44] to the Newton-Maclaurin series.

\section{Applications of the renormalization method based on Newton-Maclaurin expansion}

\subsection{A Van der Pol type of perturbation difference equation}
 Consider the equation
 \begin{equation}
 y(n+2)-2\cos\theta y(n+1)+y(n)=\epsilon(1-y^2(n+1))(y(n+2)-y(n)),
 \end{equation}
where $\epsilon$ is a positive small parameter. This equation is a center finite difference approximation of the
continuous Van der Pol differential equation[33,37]. Assume that the solution can be expanded as a power series of the
small parameter $\epsilon$
\begin{equation}
y(n)=y_0(n)+y_1(n)\epsilon+\cdots+y_k(n)\epsilon^k+\cdots.
\end{equation}
Substituting it into the above equation yields the equations of
$y_k(n)$'s such as
\begin{equation}
y_0(n+2)-2\cos\theta y_0(n+1)+y_0(n)=0,
\end{equation}
and
\begin{equation}
 y_1(n+2)-2\cos\theta y_1(n+1)+y_1(n)=(1-y_0^2(n+1))(y_0(n+2)-y_0(n)),
\end{equation}
and so forth. Solving the first equation gives
\begin{equation}
y_0(n)=A\mathrm{e}^{in\theta}+B\mathrm{e}^{-in\theta}.
\end{equation}
By the variation of constant method, we solve the second equation and obtain
\begin{equation}
y_1(n)=\{(A-A^2B)\mathrm{e}^{in\theta}+(B-AB^2)\mathrm{e}^{-in\theta}\}(n-m).
\end{equation}
Therefore, the solution is given by
\begin{equation}
y(n)=A\mathrm{e}^{in\theta}+B\mathrm{e}^{-in\theta}+
\epsilon\{(A-A^2B)\mathrm{e}^{in\theta}+(B-AB^2)\mathrm{e}^{-in\theta}\}(n-m)+O(\epsilon^2),
\end{equation}
which includes secular terms.

By the renormalization method, considering two constants $A$ and $B$ as the functions of variable $m$, and expanding $e^{in\theta}$ and $e^{-in\theta}$ as the Newton-Maclaurin series at the point $m$, a closed  renormalization equations system can be taken as
\begin{equation}
\Delta A(m)=\epsilon A(m),
\end{equation}
\begin{equation}
\Delta B(m)=\epsilon B(m).
\end{equation}
Solving the renormalization equations, we have
\begin{equation}
A(m)=A_0(1+\epsilon)^m\sim A_0\mathrm{e}^{m\epsilon},
\end{equation}
\begin{equation}
B(m)=B_0(1+\epsilon)^m\sim B_0\mathrm{e}^{m\epsilon},
\end{equation}
where $A_0$ and $B_0$ are two arbitrary constants, and in general, we take $A_0=\overline{B_0}$ to give a real solution. Then we give the global asymptotic solution
\begin{equation}
y(n)=A_0\mathrm{e}^{n(\epsilon+i\theta)}+B_0\mathrm{e}^{n(\epsilon-i\theta)}+O(\epsilon^2),
\end{equation}
which is a valid asymptotic solution for $n\sim\frac{1}{\epsilon}$.

\textbf{Remark 4}. If we take the renormalization equations as
\begin{equation}
\Delta A(m)=\epsilon A(m)(1-A(m)B(m)),
\end{equation}
\begin{equation}
\Delta B(m)=\epsilon B(m)(1-A(m)B(m)),
\end{equation}
and let $A=kB$, we have
\begin{equation}
\Delta B(m)=\epsilon B(m)(1-kB^2(m)),
\end{equation}
which has three equilibrium points $0$ and $\pm\frac{1}{\sqrt k}$, where the first is unstable and the last two are stable.
However, it is difficult to give the exact solution to the nonlinear difference equation (116)(for the similar discussions on the problem, see [37]) .

\subsection{A boundary layer problem}

Consider a singular perturbation problem[27]
\begin{equation}
\epsilon y(n+2)+a y(n+1)+by(n)=0,
\end{equation}
with the boundary conditions
\begin{equation}
y(0)=\alpha, y(N)=\beta.
\end{equation}
This is also a boundary layer problem. Assume that the solution can be expanded as a power series of the
small parameter $\epsilon$
\begin{equation}
y(n)=y_0(n)+y_1(n)\epsilon+\cdots+y_k(n)\epsilon^k+\cdots.
\end{equation}
Substituting it into the above equation yields the equations of
$y_k(n)$'s such as
\begin{equation}
a y_0(n+1)+by_0(n)=0,
\end{equation}
and
\begin{equation}
a y_1(n+1)+by_1(n)=-y_0(n+2),
\end{equation}
and so forth. Solving the first equation gives
\begin{equation}
y_0(n)=A(-\frac{b}{a})^n.
\end{equation}
The second equation becomes
\begin{equation}
ay_1(n+1)+by_1(n)=-A(-\frac{b}{a})^{n+2}.
\end{equation}
By the variation of constant method, we take
\begin{equation}
y_1(n)=B(n)(-\frac{b}{a})^n,
\end{equation}
and give
\begin{equation}
\Delta B(n)=\frac{Ab}{a^2}.
\end{equation}
Solving the above equation yields
\begin{equation}
 B(n)=B_0+\frac{Ab}{a^2}(n-m).
\end{equation}
Therefore, we have
\begin{equation}
y_1(n)=B_0(-\frac{b}{a})^n+\frac{Ab}{a^2}(n-m)(-\frac{b}{a})^n.
\end{equation}
Furthermore, the solution is given by
\begin{equation}
y(n)=A(-\frac{b}{a})^n+\epsilon B_0(-\frac{b}{a})^n+\epsilon \frac{Ab}{a^2}(n-m)(-\frac{b}{a})^n+O(\epsilon^2),
\end{equation}
which includes a secular term.

By the renormalization method, considering two constants $A$ and $B_0$ as the functions of variable $m$,  and expanding $(-\frac{b}{a})^n$  as the Newton-Maclaurin series at the point $m$, the  renormalization equation can be taken as
\begin{equation}
\Delta( A(m)+\epsilon B_0(m))=\epsilon \frac{A(m)b}{a^2}.
\end{equation}
Thus a closed renormalization equations system can be taken as
\begin{equation}
\Delta A(m)=\epsilon\frac{A(m)b}{a^2},
\end{equation}
\begin{equation}
\Delta B_0(m)=0.
\end{equation}
Solving the renormalization equations, we get
\begin{equation}
A(m)=A_0(1+\epsilon\frac{b}{a^2})^m,
\end{equation}
and $B_0$ is a constant, where $A_0$ is also an arbitrary constant. Then we give the global asymptotic solution
\begin{equation}
y(n)=A_0(1+\epsilon\frac{b}{a^2})^n(-\frac{b}{a})^n+\epsilon B_0(-\frac{b}{a})^n.
\end{equation}
According to the boundary conditions, we have
\begin{equation}
A_0+\epsilon B_0=\alpha,
\end{equation}
\begin{equation}
A_0(1+\epsilon\frac{b}{a^2})^N(-\frac{b}{a})^N+\epsilon B_0(-\frac{b}{a})^N=\beta.
\end{equation}
Solving the two equations, we obtain
\begin{equation}
A_0=\frac{\beta(-\frac{a}{b})^N-\alpha}{(1+\epsilon\frac{b}{a^2})^N-1},
\end{equation}
\begin{equation}
B_0=\frac{1}{\epsilon}\frac{\alpha(1+\epsilon\frac{b}{a^2})^N-\beta(-\frac{a}{b})^N}{(1+\epsilon\frac{b}{a^2})^N-1}.
\end{equation}

\subsection{Reduction and invariant manifold}

Here we follow [39] for notations. Consider the discrete dynamic system
\begin{equation}
\Delta X(n)=F_0(X)+\epsilon P(X,n),
\end{equation}
where $X(n)$ is a $m-$dimensional vector, and $0<\epsilon\ll 1$.  We assume that the initial point is on the invariant manifold $M$ if the invariant  manifold exists.  Further assume that $M$ is represented by the parameters vector $s$ and the reduction system of the above equation is given by the vector field $G$ as
\begin{equation}
\Delta s(n)=G(s)
\end{equation}
and the manifold $M$ is parameterized by
\begin{equation}
X(n)=R(s).
\end{equation}
Therefore our aim is to find the vector field $G$ and the function $R$. By our renormalization method based on the Newton-Maclaurin expansion, we can obtain them. Letting $X(n)=X_0(n)+\epsilon X_1(n)+\epsilon^2X_2(n)+\cdots$ and substituting it into the equation (138), we have
\begin{equation}
\Delta X_0(n)=F(X_0(n)),
\end{equation}
\begin{equation}
\Delta X_1(n)=F'(X_0(n))X_1(n)+P(X_0(n),n),
\end{equation}
and so forth. Solving the equation (141) yields
\begin{equation}
 X_0(n)=R(c_1,\cdots,c_k,n),
\end{equation}
where $c_1,\cdots,c_k$ are $k$ arbitrary constants and here $k\leq m$. Denote $C=(c_1,\cdots,c_k)$. Then the invariant manifold $M_0$ can be represented by
\begin{equation}
s(n)=C(n),
\end{equation}
and hence the reduction system will be given by the evolution equation of $C$ ( by taking $C$ as the functions of a general point $n_0$) which is just from the renormalization equation.

Now we give the method of finding the reduction equation by the following  simple example which is a discrete version of Kunihiro's example,
\begin{equation}
\Delta x(n)=\epsilon f(x(n),y(n)),
\end{equation}
\begin{equation}
\Delta y(n)=-y(n)+g(x(n)),
\end{equation}
where $f$ and $g$ are two known functions. Letting $x(n)=x_0(n)+x_1(n)\epsilon+\cdots$ and $y(n)=y_0(n)+y_1(n)\epsilon+\cdots$ and substituting them into (145) and (146),  we have
\begin{equation}
\Delta x_0(n)=0,
\end{equation}
\begin{equation}
\Delta y_0(n)=-y_0(n)+g(x_0(n)),
\end{equation}
and
\begin{equation}
\Delta x_1(n)=f(x_0(n),y_0(n)),
\end{equation}
\begin{equation}
\Delta y_1(n)=-y_1(n)+ g'(x_0(n))x_1(n).
\end{equation}
Solving the above equations give
\begin{equation}
x_0(n)=c, y_0(n)=g(c),
\end{equation}
\begin{equation}
 x_1(n)=f(c,g(c))(n-n_0)+b,
\end{equation}
\begin{equation}
 y_1(n)=g'(c)f(c,g(c))(n-n_0)+bg'(c)+g'(c)f(c,g(c)),
\end{equation}
where $c$ and $b$ are tow arbitrary constants. By our renormalization method, we assume that $c$ and $b$ are the functions of $n_0$. Then the renormalization equations are given as
\begin{equation}
 \Delta(c+\epsilon b)=\epsilon f(c,g(c)),
\end{equation}
\begin{equation}
 \Delta(g(c)+g'(c)f(c,g(c))+\epsilon bg'(c))=\epsilon g'(c)f(c,g(c)).
\end{equation}
A closed approximate renormalization equations system is taken as
\begin{equation}
 \Delta c=\epsilon f(c,g(c)),
\end{equation}
\begin{equation}
 \Delta b=0.
\end{equation}
Thus the reduction equation is just (156) and the invariant manifold is represented by
\begin{equation}
y(x)=g(x)-g'(x)f(x,g(x)).
\end{equation}

In conclusion, we have obtained the invariant manifold and the corresponding reduction dynamics of slow variables on it.

\section{Homotopy renormalization method for difference equations and applications}

The renormalization methods is not omnipotent, and has some weaknesses[1]. In order to overcome these weaknesses, we propose an iteration method as done in [1]. In fact,  what we need is only to have a freedom to
  choose the first approximate solution. The iteration method can
  give us this freedom and then can be applied to many
  equations. In the following, we will obtain this
  homotopy renormalization  method for difference equations and further give its some applications.

 Consider the difference equation
\begin{equation}
N(y(n),\Delta y(n),\cdots)=0,
\end{equation}
where $N$ is an operator. Take a simple linear equation
\begin{equation}
L(y(n),\Delta y(n),\cdots)=0,
\end{equation}
where $L$ is in general a linear operator with constant coefficients or variable coefficients. Next we take the homotopy
equation
\begin{equation}
(1-\epsilon) L(y(n),\Delta y(n),\cdots)+\epsilon N(y(n),\Delta y(n),\cdots)=0,
\end{equation}
where the homotopy parameter $\epsilon$ satisfies $0\leq \epsilon
\leq 1$. We see that the homotopy equation changes from the
simple equation $L(y(n),\Delta y(n),\cdots)=0$ to the aim equation $N(y(n),\Delta y(n),\cdots)=0$ as
$\epsilon$ changing from $0$ to $1$. Therefore, we can deform the
solution of the simple equation to that of the aim equation. Expand the solution of the homotopy equation as a power series of
$\epsilon$
\begin{equation}
y(n,\epsilon)=\sum_{k=0}^{+\infty}y_k(n)\epsilon^k,
\end{equation}
where $y_k(n)'s$ are unknown functions. Substituting the solution
into the homotopy equation and equating the coefficients of the
power of $\epsilon$, we get the linear equations
\begin{equation}
 L(y_{k+1}(n),\cdots)=L(y_k(n),\cdots)-N(y_k(n),\cdots),
\end{equation}
for $k=0,1,\cdots$. Solving these linear equations give solutions
$y_0(n),y_1(n)$ and so on.  Expanding
every $y_k(n)$ as a Newton-Maclaurin series at a general point $m$ and rearranging
the solution, and taking $\epsilon=1$ at the end, we get the solution of
original equation
\begin{equation}
y(n,m)=Y_0(m,A,B,\cdots)+Y_1(m,A,B,\cdots)(n-m)+Y_2(m,A,B,\cdots)\binom{n-m}{2}+O(1),
\end{equation}
where $A$ and $B$ and so on are parameters. From this expression, the solution is
\begin{equation}
y(m)=Y_0(m,A,B,\cdots).
\end{equation}
In order to use the renormalization  method, we  consider these
parameters as the functions of $m$ and determine them by another
relation, that is, the renormalization equation
\begin{equation}
\Delta Y_0=Y_1.
\end{equation}
Through a suitable choice of the closed equations system, we
solve out these  constants  $A$ and $B$, and
substitute them into the solution $Y_0$ to get the
asymptotic global solution.

We call the above proposed method the homotopy renormalization method based on the Newton-Maclaurin expansion, for simplicity, also the HTR method.

 In next two examples, we will see that the weaknesses of renormalization method can be overcome by the homotopy renormalization method.

\textbf{Example 1}. Consider the nonlinear difference equation
 \begin{equation}
\Delta y(n)=\eta(y(n)+y^3(n)),
\end{equation}
which cannot be dealt with by the TR method to improve the asymptotic solution by the similar reason in [1], where $\eta$ is a parameter. Here we find its nontrivial asymptotic solution by the homotopy renormalization method.

  The homotopy
equation is given by
 \begin{equation}
\Delta y(n)+\frac{1}{2}y(n)=\epsilon\{\frac{1}{2}y(n)+\eta(y(n)+y^3(n))\}.
\end{equation}
Taking
\begin{equation}
y(n)=\sum_{k=0}^{+\infty}y_k(n)\epsilon^k,
\end{equation}
and substituting it into the above homotopy equation yields
\begin{equation}
\Delta y_0(n)+\frac{1}{2}y_0(n)=0,
\end{equation}
and
\begin{equation}
\Delta y_1(n)+\frac{1}{2}y_1(n)=\frac{1}{2}y_0(n)+\eta(y_0(n)+y_0^3(n)),
\end{equation}
and so forth.

Solving the first two equations give
\begin{equation}
y_0(n)=(\frac{1}{2})^nK_0,
\end{equation}
\begin{equation}
y_1(n)=(\frac{1}{2})^m((1+2\eta)K_0+2\eta K_0^3(\frac{1}{2})^{2m})(n-m).
\end{equation}
Therefore, by takeing $\epsilon=1$, we get
\begin{equation}
y(n)=(\frac{1}{2})^nK_0+(\frac{1}{2})^m((1+2\eta)K_0+2\eta K_0^3(\frac{1}{2})^{2m})(n-m)+O(1).
\end{equation}
 According to the renormalization method, by considering $K_0$ as a function of $m$, the
renormalization equation is
\begin{equation}
\Delta K_0(m)=(1+2\eta)K_0(m)+2\eta K_0^3(m)(\frac{1}{2})^{2m},
\end{equation}
which is a nonlinear equation. For further simplicity, we take a closed approximate renormalization equation
\begin{equation}
\Delta K_0(m)=(1+2\eta)K_0(m),
\end{equation}
whose solution is
\begin{equation}
K_0(m)=2^m(1+\eta)^mB_0,
\end{equation}
where $B_0$ is a constant. Then, correspondingly, the
global solution is
\begin{equation}
y(n)=(1+\eta)^nB_0,
\end{equation}
which is uniform valid for $n\sim\frac{1}{\eta}$.

\textbf{Example 2}. Consider the second order nonlinear difference equation
\begin{equation}
y(n+2)-2y(n+1)+y(n)=D(y(n)-y^3(n)),
\end{equation}
with boundary conditions
\begin{equation}
y(0)=1, y(+\infty)=0.
\end{equation}
 This difference equation is a simplified discrete version of domain boundaries problem in convection patterns[45]. It is obvious that the equation cannot be dealt with by the TR method since $D$ is not a small parameter (in fact, we have $D=1$). Here we solve it by the HTR method.

The homotopy
equation is taken by
\begin{equation*}
y(n+1)-\frac{1+\mathrm{e}^{\lambda n}}{1+\mathrm{e}^{\lambda (n+1)}}y(n)=\epsilon\{y(n+1)-\frac{1+\mathrm{e}^{\lambda n}}{1+\mathrm{e}^{\lambda (n+1)}}y(n)
\end{equation*}
\begin{equation}
+k(y(n+2)-2y(n+1)+y(n)-D(y(n)-y^3(n)))\}.
\end{equation}
Taking
\begin{equation}
y(n)=\sum_{k=0}^{+\infty}y_k(n)\epsilon^k,
\end{equation}
and substituting it into the above homotopy equation yields
\begin{equation}
y_0(n+1)-\frac{1+\mathrm{e}^{\lambda n}}{1+\mathrm{e}^{\lambda (n+1)}}y_0(n)=0,
\end{equation}
and
\begin{equation}
y_1(n+1)-\frac{1+\mathrm{e}^{\lambda n}}{1+\mathrm{e}^{\lambda (n+1)}}y_1(n)=k(y_0(n+2)-2y_0(n+1)+y_0(n)-D(y_0(n)-y_0^3(n)),
\end{equation}
and so forth.

Solving the two first equations give
\begin{equation}
y_0(n)=\frac{A}{1+\mathrm{e}^{\lambda n}},
\end{equation}
\begin{equation}
y_1(n)=f(m)(n-m)+\cdots,
\end{equation}
where
\begin{equation}
f(n)=k(y_0(n+2)-2y_0(n+1)+y_0(n)-D(y_0(n)-y_0^3(n))).
\end{equation}
Therefore, by taking $\epsilon=1$, we get
\begin{equation}
y(n)=\frac{A}{1+\mathrm{e}^{\lambda n}}+f(m)(n-m)+O(1).
\end{equation}
 According to the renormalization method, by considering $A$ as a function of $m$, a closed
renormalization equation can be taken as
\begin{equation}
\Delta A(m)=k(1-D)A(m),
\end{equation}
whose solution is
\begin{equation}
A(m)=A_0(k(1-D)+1)^m,
\end{equation}
where $A_0$ is a constant. Then, correspondingly, the
global approximate asymptotic solution is
\begin{equation}
y(n)=\frac{A_0(k(1-D)+1)^n}{1+\mathrm{e}^{\lambda n}}.
\end{equation}
According to the boundary conditions, we take
\begin{equation}
A_0=2,D=1, \lambda>0.
\end{equation}
Therefore, the solution becomes
\begin{equation}
y(n)=\frac{2}{1+\mathrm{e}^{\lambda n}},
\end{equation}
which is a global valid asymptotic solution satisfying boundary conditions.

\section{Conclusions}

The renormalization method  based on the Newton-Maclaurin series is proposed to give the asymptotic solutions to difference equations such as  regular and singular perturbation equations. Furthermore, to overcome some weaknesses of the renormalization method, a homotopy renormalization method is introduced to deal with some difference equations including these problems without small parameters. As a result, we can see that the  renormalization method is clear in theory and simple in practice. It is undoubt that the proposed remormalization methods can also be applied to more linear and nonlinear difference and differential equations such as  systems with delay[46], pantograph equations[47], lattice equations[48],  the difference equations systems and partial difference equations.

\end{document}